\documentclass{SCIS2019}
\usepackage{CJKutf8}
\usepackage{setspace}
\usepackage{mathrsfs}
\usepackage{graphicx,enumerate,cleveref}
\usepackage{accents,amssymb,amsmath,wasysym,cite,color,float,bm}
\newcommand{\supercite}[1]{\textsuperscript{\cite{#1}}}
\makeatletter
\def\widebar{\accentset{{\cc@style\underline{\mskip10mu}}}}
\makeatother

\usepackage{geometry}

\usepackage{setspace}

\allowdisplaybreaks

\usepackage{paralist}

\makeatletter
\renewcommand*\env@matrix[1][\arraystretch]{%
	\edef\arraystretch{#1}%
	\hskip -\arraycolsep
	\let\@ifnextchar\new@ifnextchar
	\array{*\c@MaxMatrixCols c}}
\makeatother

\usepackage{fancyhdr} 
\usepackage{bm}

\begin{document}
\begin{CJK}{UTF8}{gbsn}
\ArticleType{NOTE}
\Year{2019}
\Month{}
\Vol{}
\No{}
\DOI{}
\ArtNo{}
\ReceiveDate{}
\ReviseDate{}
\AcceptDate{}
\OnlineDate{}

\setlength{\footskip}{20pt}

\title{A Note on Global Optimization for Max-Plus Linear Systems}{A note on global optimization for max-plus linear systems}

\author[1]{Cailu WANG}{}
\author[2]{Yuegang TAO}{{yuegangtao@hebut.edu.cn}}


\AuthorCitation{Wang C and Tao Y}


\address[1]{School of Automation, Beijing Institute of Technology, Beijing {\rm 100081}, China}
\address[2]{School of Artificial Intelligence, Hebei University of Technology, Tianjin {\rm 300130}, China}

\abstract{This note further addresses the global optimization problem for max-plus linear systems considered in [Automatica 119 (2020) 109104]. Firstly, 
	the operations between $\pm\infty$ and real numbers involved in the formulas of solving global optimization problems are explained explicitly. Secondly, the formula of the greatest lower bound and the criterion of solvability  of globally optimal solutions are simplified.
	Thirdly, the criterion of uniqueness of globally optimal solutions and the set of all globally optimal solutions are presented.}


\maketitle

\begin{multicols}{2}

\section{Introduction}\vspace{-3mm}
Max-plus linear systems can describe
some nonlinear time-evolution systems with synchronization but no concurrency, such as flexible manufacturing systems, flow shop scheduling, traffic managements,
communication networks (see, e.g., \cite{8,MaxPlusAtWork,SynchronizationLinearity}).
Many progresses have been made in 
control and optimization of max-plus linear systems (see, e.g., \cite{Farahani2017On,PMPG,Optimistic,MYSKOVA20208,Goncalves2019on,screening,IJRNC,LinearSystemTheoretic,8993726,Candido2018Conditional,PeriodicManufacturing,ModelPredictive,Ioannis2018Stochastic,Kubo2018Applications,Mendes2019Stochastic,7335578}).

Global optimization is to find the global minimizer of a
function or a set of functions over a given set, which has been
a basic tool in all areas of engineering, medicine, economics,
and other sciences. The 
global optimization for max-plus linear systems has been considered in \cite{Tao2020Global}, whose objective function is a max-plus function and constraint function is a real function. 
This note is a further explanation and extension of the results in \cite{Tao2020Global}.

Some formulas for globally optimal solutions in\\ \cite{Tao2020Global} involves the operations between some special elements in max-plus algebra and real numbers. 
This note will explain these operations definitely and simplify the formulas. 
In addition, reference \cite{Tao2020Global} neglects to consider the $0$ coefficients in the
\begin{spacing}{0.1}
{\color{white}1}
\end{spacing} \noindent 
  constraint function when it discusses the uniqueness of globally optimal solutions. 
This note will present a necessary and sufficient condition for the uniqueness and construct the set of all globally optimal solutions.




\vspace{-4mm}
\section{Preliminaries} \vspace{-2mm}
Let us introduce some basic definitions and notations from the max-plus algebra, which can
be consulted for more details in \cite{8,MaxPlusAtWork,SynchronizationLinearity}.

Let $\mathbb{R}$ be the set of real numbers, $\mathbb{N}$ be the set of natural numbers and $\mathbb{N}^{+}$ be the set of positive integers. {For $n\in\mathbb{N}^{+}$,
	denote by $\mathbb{N}_n$ the set $\{1, 2,\ldots,$ $ n\}$.}
Let $\mathbb{R}_{\max}$ be the set $\mathbb{R} \cup\{-\infty\}$ with max and $\mathbb{+}$ as two binary operations $\oplus$ and $\otimes$, respectively, i.e., for $a,b\in\mathbb{R}_{\max }$,
$$\setlength{\abovedisplayskip}{5pt}
\setlength{\belowdisplayskip}{5pt}
a\oplus b=\max\{a,b\}\text{\ and\ } a\otimes b=a+b.$$
$\left(\mathbb{R}_{\max }, \oplus, \otimes\right)$ is called the {\it max-plus algebra}, in whi- ch $-\infty$ is the zero element denoted by $\varepsilon$, and 0 is the identity element denoted by $e$. 
For $a,b\in\mathbb{R}_{\max}$, $a\leqslant b$ if $a\oplus b=b$.

Let $\mathbb{R}_{\max}^{n}$ and $\mathbb {R}_{\max }^{m \times n}$ be the set of $n$-dimensional vectors and $m \times n$ matrices with entries in $\mathbb{R}_{\mathrm{max}}$, respectively. 
The vectors and matrices are represented by bold-type letters.
The addition and multiplication
can be extended to the vector-matrix algebra over $\mathbb{R}_{\max}$ by
the direct analogy to the con- ventional linear algebra: For $\bm{A}=(a_{i j}),\ \bm{B}=(b_{i j})$ $\in \mathbb{R}_{\max }^{m \times n}$,
$(\boldsymbol{A} \oplus \boldsymbol{B})_{i j}=a_{i j} \oplus b_{i j};$
for $\bm{A}=(a_{i j}) \in$ $ \mathbb{R}_{\max }^{m \times p}$ and 
$\boldsymbol{B}=(b_{i j}) \in \mathbb{R}_{\max }^{p \times n}$, 
\[\setlength{\abovedisplayskip}{5pt}
\setlength{\belowdisplayskip}{5pt}
(\boldsymbol{A} \otimes \boldsymbol{B})_{i j}=\bigoplus_{k=1}^{p} a_{i k} \otimes b_{k j}.\]
For 
$\bm{A},\bm{B}$ $\in\mathbb{R}_{\max}^{m\times n}$,
$\bm{A}\leqslant \bm{B}$ if $\bm{A}\oplus \bm{B}=\bm{B}$.
For $\bm{A},$ $\bm{B}\in\mathbb{R}_{\max}^{m\times n}$ and $\bm{x}, \bm{y}\in\mathbb{R}_{\max}^{n}$,
if $\bm{A}\leqslant \bm{B}$ and $\bm{x}\leqslant \bm{y}$, then $\bm{A}\otimes\bm{x}$ $\leqslant \bm{B}\otimes\bm{y}.$


For $\bm{A}\in\mathbb{R}_{\max}^{m\times n}$ and $\bm{x}\in\mathbb{R}_{\max}^{n}$,
$F(\bm{x})=\bm{A}\otimes \bm{x}$ is called a {\it max-plus function} of type $(n, m)$.
Denoted by $F_{i}(\bm{x})$
the $i$th component of $F(\bm{x})$.
A {\it max-plus linear system} is a system that can be described by using max-plus functions.

Given $\bm{A}=(a_{ij})\in\mathbb{R}_{\max}^{m\times n}$ and  
$\bm{b}=(b_{i})\in\mathbb{R}_{\max}^{m}$, 
a \textit{system of max-plus linear equations} is defined as
\begin{equation}\label{e1}
\setlength{\abovedisplayskip}{5pt}
\setlength{\belowdisplayskip}{5pt}
\bm{A}\otimes \bm{x}=\bm{b},
\end{equation}
where $\bm{x}=(x_{j})$ is considered in $\widebar{\mathbb{R}}_{\max} (=\mathbb{R}_{\max}\cup\{+\infty\})$, rather than in $\mathbb{R}_{\max}$. 
System (\ref{e1}) is said to be {\it solvable}
if there exists $\tilde{\bm{x}}\in$ $\widebar{\mathbb{R}}_{\max}^{n}$ such that $\bm{A}\otimes\tilde{\bm{x}}=\bm{b}$, and $\tilde{\bm{x}}$ is called a \textit{solution} of system (\ref{e1}). 
For $\tilde{\bm{x}}\in\widebar{\mathbb{R}}_{\max}^{n}$, $\tilde{\bm{x}}$ is called a \textit{subsolution} of system (\ref{e1}) if  $\bm{A}\otimes\tilde{\bm{x}}\leqslant\bm{b}$.

\vspace{1.5mm}
\begin{lemma}\label{l1}\supercite{8}
	For system {(\ref{e1})}, let $\bm{x}^{*}(\bm{A},\bm{b})\in\widebar{\mathbb{R}}_{\max}^{n}$ be defined by
	\begin{equation}\label{e10}
	\setlength{\abovedisplayskip}{5pt}
	\setlength{\belowdisplayskip}{5pt}
	x_{j}^{*}(\bm{A},\bm{b})=\min\limits_{i\in\mathbb{N}_{m}}\{b_{i}-a_{ij}\},\ j\in\mathbb{N}_{n}.
	\end{equation}
	Then, $\bm{x}^{*}(\bm{A},\bm{b})$ is the {greatest subsolution} of system {(\ref{e1})}, i.e.,
	\begin{enumerate}[i)\leftmargin=-10pt]
		\item $\bm{A}\otimes \bm{x}^{*}(\bm{A},\bm{b})\leqslant \bm{b}$;
		\item if $\tilde{\bm{x}}$ is a subsolution of system {(\ref{e1})}, then $\tilde{\bm{x}}\leqslant \bm{x}^{*}(\bm{A},\bm{b})$. In particular, if $\tilde{\bm{x}}$ is a solution,
		then $\tilde{\bm{x}}\leqslant \bm{x}^{*}(\bm{A},\bm{b})$. 
	\end{enumerate} 
\end{lemma}\vspace{1.5mm}

System {(\ref{e1})} is solvable if and only if the greatest subsolution is a solution, i.e., $\bm{A}\otimes \bm{x}^{*}(\bm{A},\bm{b})= \bm{b}$.

\vspace{-3mm}
\section{Global optimization}\vspace{-2mm}
Consider the \textit{global optimization problem}
\begin{equation}\label{e2}
\setlength{\abovedisplayskip}{5pt}
\setlength{\belowdisplayskip}{5pt}
\min_{\bm{x}\in\mathcal{X}}  F(\bm{x}),
\end{equation}
where the optimization variable is $\bm{x}=(x_{j})\in\mathbb{R}^{n}$;
the objective function is $F(\bm{x})=\bm{A}\otimes \bm{x}$, in which  $\bm{A}=(a_{ij})\in\mathbb{R}_{\max}^{m\times n}$;
the constraint set is
\[\setlength{\abovedisplayskip}{5pt}
\setlength{\belowdisplayskip}{5pt}	\mathcal{X} =\left\{\bm{x}\,\Big|\displaystyle\sum\limits_{j=1}^{n}{k_{j} x_{j}}=c,\ k_{j}\geqslant 0,\ 
k_{j}\text{\,are not all 0}\right\},\]
in which $c\in\mathbb{R}$ is a constant. 

\vspace{1.5mm}
\begin{definition}\label{d1}\supercite{Tao2020Global}
	Problem \eqref{e2} is said to be {\it solvable}
	if there exists $\tilde{\bm{x}}\in\mathcal{X}$ such that $F(\bm{x})\geqslant F(\tilde{\bm{x}})$ for any $\bm{x}$ $\in\mathcal{X}$, and
	$\tilde{\bm{x}}$ is called a {\it globally optimal solution} of problem \eqref{e2}. Otherwise, problem \eqref{e2} is said to be {\it unsolvable}.
\end{definition}\vspace{1.5mm}

Note that, for problem \eqref{e2}, an $\bm{\varepsilon}$ row in $\bm{A}$ has no affect on globally optimal solutions. 
In addition, if the $j_{0}$th column of $\bm{A}$ is $\bm{\varepsilon}$ and meanwhile $k_{j_{0}}=0$, 
then $x_{j_{0}}$ disappears from both the objective function and the constraint function of problem \eqref{e2}. 
Consequently, problem \eqref{e2} becomes into the optimization of variables $\{x_{j}\,|\,j\neq j_{0}\}$. 
Without loss of generality, let us make the following assumptions:
\begin{enumerate}[i)\leftmargin=-10pt]
	\item $\bm{A}$ is row $\mathbb{R}$-astic, i.e., $\bm{A}$ has no $\bm{\varepsilon}$ row.
	\item $k_{j}$$\neq$0 if the $j$th column of $\bm{A}$ is $\bm{\varepsilon}$.
\end{enumerate}

\vspace{1.5mm}
\begin{theorem}\label{l2}\supercite{Tao2020Global}
	For problem \eqref{e2}, let $\bm{b}=(b_i)\in\mathbb{R}_{\max}^{m}$ be defined by
	\begin{equation}\label{e04}
	\setlength{\abovedisplayskip}{5pt}
	\setlength{\belowdisplayskip}{5pt}
	b_{i}=\frac{\sum_{j=1}^{n}k_{j}a_{ij}+c}{\sum_{j=1}^{n}k_{j}},\ i\in{\mathbb{N}_{m}}.
	\end{equation}
	Then, $\bm{b}$ is the greatest lower bound of problem \eqref{e2}:
	\begin{enumerate}[i)\leftmargin=-10pt]
		\item $F(\bm{x})\geqslant \bm{b}$ for any $\bm{x}\in\mathcal{X}$;
		\item if $\bm{y}$ is a lower bound of problem \eqref{e2}, then $\bm{y}\leqslant \bm{b}$.
	\end{enumerate}
\end{theorem}\vspace{1.5mm}

The calculation of \eqref{e04} may refer to some opera- tions between elements $\pm \infty$  in $\widebar{\mathbb{R}}_{\max}$ and real numbers in $\mathbb{R}$. 
To compute the greatest lower bound by \eqref{e04}, 
one needs to make the following statements: For $r,p\in\mathbb{R}$ and $p>0$,  \vspace{-4mm}
\begin{center}\renewcommand\arraystretch{1.2}
	\begin{tabular}{cc}
		$\max\{r,-\infty\}=r$, & $\max\{r,+\infty\}=+\infty$,\\
		$r+(-\infty)=-\infty$,& $r+(+\infty)=+\infty$,\\ 
		$r-(-\infty)=+\infty$,& $r-(+\infty)=-\infty$,\\
		$(-\infty)-r=-\infty$,& $(+\infty)-r=+\infty$,\\
		$p\times (-\infty)=-\infty$,& $p\times (+\infty)=+\infty$,\\ 
		$0\times (-\infty)=0$,& $0\times (+\infty)=0$,\\
		$(+\infty)+(-\infty)=-\infty$,& $(-\infty)-(-\infty)=+\infty$.\\
	\end{tabular}
\end{center}

Considering that $0$ coefficients in the constraint function have no function in calculating the greatest lower bound, formula \eqref{e04} can equivalently be simplified as
\begin{equation}\label{e4}
\setlength{\abovedisplayskip}{5pt}
\setlength{\belowdisplayskip}{5pt}
b_{i}=\frac{\sum_{j\in \mathcal{J}}k_{j}a_{ij}+c}{\sum_{j\in \mathcal{J}}k_{j}},\ i\in{\mathbb{N}_{m}},
\end{equation}
where
$
\mathcal{J}=\{j\in\mathbb{N}_{n}\,|\,k_{j}>0\}.
$
The formula above avoids the multiplication $\times$ between $0$ and infinity elements in $\widebar{\mathbb{R}}_{\max}$.
\newpage
\vspace{-3mm}
\section{Optimal solutions}\vspace{-2mm}
Let us provide two necessary conditions for the existence of globally optimal solutions.

\vspace{1.5mm}
\begin{lemma}\label{cor1}
If problem \eqref{e2} is solvable, then 
	\begin{enumerate}[i)\leftmargin=-10pt]
	\item the greatest lower
	bound $\bm{b}$ given in \eqref{e4} is finite;
	\item $\bm{x}^{*}(\bm{A},\bm{b})$ is finite.
\end{enumerate}
\end{lemma}\vspace{1.5mm}\noindent
\textit{Proof.} 1) Suppose that there exists $i_{0}\in\mathbb{N}_{m}$ such that $b_{i_{0}}=\varepsilon$. 
Since $\bm{A}$ is row $\mathbb{R}$-astic, there exists $j_{0}\in$ $\mathbb{N}_{n}$ such that $a_{i_{0}j_{0}}\neq\varepsilon$. 
Let $\tilde{\bm{x}}=(\tilde{x}_{j})$ be a globally optimal solution of problem \eqref{e2}.
Then, 
$$\setlength{\abovedisplayskip}{5pt}
\setlength{\belowdisplayskip}{5pt}\varepsilon=b_{i_{0}}=F_{i_{0}}(\tilde{\bm{x}})=\max_{j\in\mathbb{N}_{n}}\{a_{i_{0}j}+\tilde{x}_{j}\}\geqslant a_{i_{0}j_{0}}+\tilde{x}_{j_{0}}.$$
This implies that $\tilde{x}_{j_{0}}=\varepsilon$, and hence $\tilde{\bm{x}}\notin\mathcal{X}$. 
This contradiction indicates that  $\bm{b}$ is finite.

2) Since $\bm{b}$ is finite,  ${x}_{j}^{*}(\bm{A},\bm{b})\neq\varepsilon$ for any $j\in\mathbb{N}_{n}$. 
Suppose that  
there exists $j_{0}\in\mathbb{N}_{n}$ such that 
$$\setlength{\abovedisplayskip}{5pt}
\setlength{\belowdisplayskip}{5pt}{x}_{j_{0}}^{*}(\bm{A},\bm{b})=\min\limits_{i\in\mathbb{N}_{m}}\{b_{i}-a_{ij_{0}}\}=+\infty.$$
This implies that $a_{ij_{0}}=\varepsilon$ for any $i\in\mathbb{N}_{m}$. 
Assum- ption 2) ensures that $k_{j_{0}}>0$, i.e, $j_{0}\in\mathcal{J}$. 
According to \eqref{e4}, $\bm{b}=\bm{\varepsilon}$, and hence problem \eqref{e2} is unsolvable. This contradiction indicates that $\bm{x}^{*}(\bm{A},\bm{b})$ is finite.
$\Square$
\vspace{1.5mm}

It can be inferred from the lemma above that if 
\begin{equation}\label{e7}
\setlength{\abovedisplayskip}{5pt}
\setlength{\belowdisplayskip}{5pt}
\sum_{j\in\mathcal{J}}k_{j}x_{j}^{*}(\bm{A},\bm{b})\big(=\sum_{j=1}^{n}k_{j}x_{j}^{*}(\bm{A},\bm{b})\big)=c,
\end{equation}
then $\bm{x}^{*}(\bm{A},\bm{b})$ is finite and $\bm{x}^{*}(\bm{A},\bm{b})\in\mathcal{X}$.
The cri-\\ teria of solvability given in \cite{Tao2020Global} can be then represented as follows.

\vspace{1.5mm}
\begin{theorem}\label{l3}\supercite{Tao2020Global}
	For problem \eqref{e2}, the following statements are equivalent:
	\begin{enumerate}[i)\leftmargin=-10pt]
	\item Problem \eqref{e2} is solvable.
	\item Equation \eqref{e7} holds.
	\item $\bm{x}^{*}(\bm{A},\bm{b})$ is a globally optimal solution. 
\end{enumerate}	
\end{theorem}\vspace{1.5mm}

Let us see a numerical example of global optimization problem that contains $\varepsilon$ coefficients in the objective function and $0$ coefficients in the constraint function.
\vspace{1.5mm}
\begin{example}\label{ex1}
	Consider the global optimization\\ problem
\begin{equation}\label{e14}
\setlength{\abovedisplayskip}{5pt}
\setlength{\belowdisplayskip}{5pt}
\min_{\bm{x}\in\mathcal{X}}  F(\bm{x}),
\end{equation}
where the objective function is $F(\bm{x})=\bm{A}\otimes \bm{x}$,
$$\setlength{\abovedisplayskip}{5pt}
\setlength{\belowdisplayskip}{5pt}\bm{A}=\begin{pmatrix}
1& 2&-2\\
-1&e&\varepsilon\\
e &1&3
\end{pmatrix},$$
and the constraint set is $\mathcal{X}=\{\bm{x}\,|\,2x_{1}+x_{2}=2\}.$
By \eqref{e4}, $\bm{b}=(2\ \,e\ \,1)^{\intercal}.$ 
By \eqref{e10}, $\bm{x}^{*}(\bm{A},\bm{b})=(1\ \,e\ -2)^{\intercal}.$
	Then,
$2x_{1}^{*}(\bm{A},\bm{b})+x_{2}^{*}(\bm{A},\bm{b})=2=c.$
It follows from \Cref{l3} that problem \eqref{e14} is solvable, and 
$(1\ \,e$ $-2)^{\intercal}$ is a globally optimal solution.

If the constraint set of problem \eqref{e14} is replaced by $\widebar{\mathcal{X}}=\{\bm{x}\,|\,2x_{1}+x_{2}+x_{3}=2\}$, then $b_{2}=\varepsilon$. It can be obtained from \Cref{cor1} that such a problem has no globally optimal solution.
\end{example}\vspace{1.5mm}


Note that, the globally optimal solution introduced in \Cref{d1} may be not unique.
According to \cite[Theorem 3]{Tao2020Global}, if problem \eqref{e2} is solvable, then
$k_{j}>0$ for any $j\in\mathbb{N}_{n}$ is a sufficient condition for the uniqueness of
globally optimal solutions. 
In fact, such a condition is also necessary.

\vspace{1.5mm}
\begin{theorem}\label{th1}
	Let problem \eqref{e2} be solvable. Problem \eqref{e2} has a unique 
	globally optimal solution if and only if $k_{j}>0$ for any $j\in\mathbb{N}_{n}$.
\end{theorem}\vspace{1.5mm}\noindent
\textit{Proof.}
	The sufficiency can be obtained from  the proof of \cite[Theorem 3]{Tao2020Global}. 
	Let us now prove the ne-\\ cessity. Since problem
	\eqref{e2} is solvable, equation \eqref{e7} holds. 
	Suppose that there exists $j_{0}\in\mathbb{N}_{n}$ such that\\ $k_{j_{0}}=0$, 
	i.e., $j_{0}\notin{\mathcal{J}}$. 
	Let $\tilde{\bm{x}}$ $=(\tilde{x}_j)\in\mathbb{R}_{\max}^{n}$ be defined by 
	\begin{equation}\label{eq2}
	\setlength{\abovedisplayskip}{5pt}
	\setlength{\belowdisplayskip}{5pt}
		\tilde{x}_{j}=\left\{\begin{array}{ll}
			\varepsilon,&j=j_{0};\\
			x^{*}_{j}(\bm{A},\bm{b}),&j\neq j_{0}.
		\end{array}\right.
	\end{equation}
It can be known from \Cref{cor1} that  $\bm{x}^{*}(\bm{A},\bm{b})$ is finite. 
	Then, $\tilde{\bm{x}}<\bm{x}^{*}(\bm{A},\bm{b})$ and 
	\begin{align*}
	\setlength{\abovedisplayskip}{5pt}
	\setlength{\belowdisplayskip}{5pt}
		\sum_{j\in\mathcal{J}}k_{j}\tilde{x}_{j}
		&=\sum_{j\in\mathcal{J}}k_{j}x^{*}_{j}(\bm{A},\bm{b})=c,
	\end{align*}
	i.e., $\tilde{\bm{x}}\in\mathcal{X}$. By \Cref{l2}, $F(\tilde{\bm{x}})\geqslant \bm{b}$. 
	In addition, 
	\begin{align*}
		\setlength{\abovedisplayskip}{5pt}
		\setlength{\belowdisplayskip}{5pt}
		F(\tilde{\bm{x}})&\leqslant F(\bm{x}^{*}(\bm{A},\bm{b}))=\bm{A}\otimes\bm{x}^{*}(\bm{A},\bm{b})\leqslant \bm{b}.
	\end{align*}
	Hence, $F(\tilde{\bm{x}})= \bm{b}$. This implies that $\tilde{\bm{x}}$ is also a globally optimal solution, which contradicts with the uniqueness. 
	Hence,  $k_{j}>0$ for any $j\in\mathbb{N}_{n}$.
$\Square$

\vspace{1.5mm}
The proof of Theorem \ref{th1} is constructive and
formula \eqref{eq2} can be extended to find all globally optimal
solutions as follows.

\vspace{1.5mm}
\begin{theorem}\label{th2}
	If problem \eqref{e2} is solvable, then the set of all globally optimal solutions is 
\begin{align*}
\setlength{\abovedisplayskip}{5pt}
\setlength{\belowdisplayskip}{5pt}
\mathcal{S}=\left\{\tilde{\bm{x}}\ \Bigg|\begin{array}{l}
\tilde{x}_{j}\leqslant x^{*}_{j}(\bm{A},\bm{b}) \text{\ if\ }k_{j}=0;\\
\tilde{x}_{j}= x^{*}_{j}(\bm{A},\bm{b}) \text{\ if\ }k_{j}> 0.
\end{array}\right\}.
\end{align*}
\end{theorem}
\noindent
\textit{Proof.}
Since problem \eqref{e2} is solvable, equation \eqref{e7} holds. 
On the one hand, for any $\tilde{\bm{x}}\in\mathcal{S}$,
\begin{align*}
\setlength{\abovedisplayskip}{5pt}
\setlength{\belowdisplayskip}{5pt}
\sum_{j\in \mathcal{J}}k_{j}\tilde{x}_{j}&=\sum_{j\in \mathcal{J}}k_{j}x^{*}_{j}(\bm{A},\bm{b})=c,
\end{align*}
i.e., $\tilde{\bm{x}}\in\mathcal{X}$. By \Cref{l2}, $F(\tilde{\bm{x}})\geqslant \bm{b}$. 
In addition, 
\begin{align*}
\setlength{\abovedisplayskip}{5pt}
\setlength{\belowdisplayskip}{5pt}
F(\tilde{\bm{x}})&\leqslant F(\bm{x}^{*}(\bm{A},\bm{b}))=\bm{A}\otimes\bm{x}^{*}(\bm{A},\bm{b})\leqslant \bm{b}.
\end{align*}
Hence, $F(\tilde{\bm{x}})= \bm{b}$. This implies that $\tilde{\bm{x}}$ is a globally optimal solution. 
On the other hand, let $\tilde{\bm{x}}$ be a globally optimal solution.  
Then, $\tilde{\bm{x}}\in\mathcal{X}$ and $F(\tilde{\bm{x}})$ $=\bm{A}\otimes\tilde{\bm{x}}=\bm{b}.$ 
By \Cref{l1}, $\tilde{\bm{x}}\leqslant\bm{x}^{*}(\bm{A},\bm{b})$.
Suppose that there exists 
$j_{0}\in{\mathcal{J}}$ such that  
$\tilde{x}_{j_{0}}<x^{*}_{j_{0}}(\bm{A},\bm{b})$. Then, \vspace{-1mm}
\begin{align*}
\setlength{\abovedisplayskip}{5pt}
\setlength{\belowdisplayskip}{5pt}
c&=\sum_{j\in\mathcal{J}\setminus\{j_{0}\}}k_{j}\tilde{x}_{j}+k_{j_{0}}\tilde{x}_{j_{0}}\\&<\sum_{j\in\mathcal{J}\setminus\{j_{0}\}}k_{j}x^{*}_{j}(\bm{A},\bm{b})+k_{j_{0}}x^{*}_{j_{0}}(\bm{A},\bm{b})\\&=\sum_{j\in \mathcal{J}}k_{j}x^{*}_{j}(\bm{A},\bm{b})=c.
\end{align*}This contradiction implies that $\tilde{x}_{j}=x^{*}_{j}(\bm{A},\bm{b})$ for any 
$j$ $\in{\mathcal{J}}$, and so $\tilde{\bm{x}}\in\mathcal{S}$. Hence, $\mathcal{S}$ is the set of all
globally optimal solutions.
$\Square$\vspace{1.5mm}

Owing to the theorem above, $\bm{x}^{*}(\bm{A},\bm{b})$ is referred to as the \textit{greatest globally optimal solution} of problem \eqref{e2} if equation \eqref{e7} holds.

\vspace{1.5mm}
\begin{example}
Find all globally optimal solutions of problem \eqref{e14}. It has been calculated in \Cref{ex1} that $\bm{x}^{*}(\bm{A},\bm{b})=(1\ \,e\ -2)^{\intercal}.$ By \Cref{th2}, the set of  all globally optimal solutions is $$\setlength{\abovedisplayskip}{5pt}
\setlength{\belowdisplayskip}{5pt}\mathcal{S}=\left\{\,\left(1\ \ 0\ \ a \right)^{\intercal} \Big|\,a\leqslant -2\right\}.$$
\end{example}


\Acknowledgements{The authors would like to thank the peer experts for their
	valuable comments and sugges- tions. } 

\end{multicols}

\end{CJK}

\begin{thebibliography}{99}
	

	
\bibitem{SynchronizationLinearity}
F.~Baccelli, G.~Cohen, G.~J. Olsder, and J.~P. Quadrat.
\newblock {\em Synchronization and Linearity: An Algebra for Discrete Event
	Systems}.
\newblock John Wiley and Sons, New York, 1992.


\bibitem{MaxPlusAtWork}
B.~Heidergott, G.~J. Olsder, and J.~van~der Woude.
{\em Max-Plus at Work: Modeling and Analysis of Synch- ronized Systems}.
Princeton University Press, New Jersey, 2006.


\bibitem{8}
P.~Butkovi$\check{\text{c}}$.
{\em Max-Linear Systems: Theory and Algorithms}.
Springer-Verlag, Berlin, 2010.



\bibitem{LinearSystemTheoretic}
G.~Cohen, D.~Dubois, J.~P. Quadrat, and M.~Viot.
A linear-system-theoretic view of discrete-event processes and its
use for performance evaluation in manufacturing.
{\em IEEE Transactions on Automatic Control}, {30}:\ 210--220,
1985.

\bibitem{ModelPredictive}
B.~D{e Schutter} and T.~van~den Boom.
Model predictive control for max-plus-linear discrete event systems.
{\em Automatica}, {37}(7):\ 1049--1056, 2001.




\bibitem{PMPG}
S.~Gaubert, R.~D. Katz, and S.~Sergeev.
Tropical linear-fractional programming and parametric mean payoff
games.
{\em Journal of Symbolic Computation}, {47}:\ 1447--1478,
2012.



\bibitem{screening}
Y.~Shang, L.~Hardouin, M.~Lhommeau, and C.~A. Maia.
An integrated control strategy to solve the disturbance decoupling
problem for max-plus linear systems with applications to a high throughput
screening system.
{\em Automatica}, {63}:\ 338--348, 2016.

\bibitem{Optimistic}
J.~Xu, T.~van~den Boom, and B.~De Schutter.
Optimistic optimization for model predictive control of max-plus
linear systems.
{\em Automatica}, {74}:\,16--22, 2016.

\bibitem{7335578}
S.~E.~Z. Soudjani, D.~Adzkiya, and A.~Abate.
Formal verification of stochastic max-plus-linear systems.
{\em IEEE Transactions on Automatic Control}, {61}(10):\
2861--2876, 2016.

\bibitem{Farahani2017On}
S.~S. Farahani, T.~van~den Boom, and B.~{De Schutter}.
On optimization of stochastic max-min-plus-scaling systems:\,An
approximation approach.
{\em Automatica}, {83}: 20--27, 2017.

	\bibitem{PeriodicManufacturing}
B.~Cottenceau, L.~Hardouin, and J.~Trunk.
\newblock Weight-balanced timed event graphs to model periodic phenomena in
manufacturing systems.
\newblock {\em IEEE Transactions on Automation Science and Engineering}, 
{14}:\ 1731--1742, 2017.

	\bibitem{Candido2018Conditional}
R.~M.~F. Candido, L.~Hardouin, M.~Lhommeau, R.~S. Mendes.
Conditional reachability of uncertain max plus linear systems.
{\em Automatica}, {94}:\ 426--435, 2018.

	\bibitem{Ioannis2018Stochastic}
K.~Ioannis, P.~Maragos, and G.~Papavassilopoulos.
Stochastic stability in max-product and max-plus systems with
markovian jumps.
{\em Automatica}, {92}:\ 123--132, 2018.

\bibitem{Kubo2018Applications}
S.~Kubo and K.~Nishinari.
Applications of max-plus algebra to flow shop scheduling problems.
{\em Discrete Applied Mathematics}, {247}:\ 278--293, 2018.

	\bibitem{IJRNC}
C.~Wang, Y.~Tao, and H.~Yan.
\newblock Optimal input design for uncertain max-plus linear systems.
\newblock {\em International Journal of Robust and Nonlinear Control}, {\bf
	28}:\ 4816--4830, 2018.

\bibitem{Goncalves2019on}
V.~M. Goncalves, C.~A. Maia, and L.~Hardouin.
On max-plus linear dynamical system theory: The observation problem.
{\em Automatica}, {107}:\ 103--111, 2019.

	\bibitem{Mendes2019Stochastic}
R.~S. Mendes, L.~Hardouin, and M.~Lhommeau.
Stochastic filtering of max-plus linear systems with bounded
disturbances.
{\em IEEE Transactions on Automatic Control}, {64}(9):\
3706--3715, 2019.


\bibitem{MYSKOVA20208}
H. My\v{s}kov\'a and J. Plavka.
Interval robustness of (interval) max-plus matrices.
\emph{Discrete Applied Mathematics}, 284: 8--19, 2020.

\bibitem{8993726}
E.~Berthier and F.~Bach.
\newblock Max-plus linear approximations for deterministic continuous-state
markov decision processes.
\newblock {\em IEEE Control Systems Letters}, {4}:\ 767-- 772, 2020.


\bibitem{Tao2020Global} 
Y. Tao and C. Wang. 
Global optimization for max-plus linear systems and applications in distributed systems.
{\em Automatica}, {119}:\ 109104, 2020.



\end{thebibliography}
\end{document}